# Homological mirror symmetry for Brieskorn-Pham singularities

Masahiro Futaki and Kazushi Ueda

**Abstract**

We prove that the derived Fukaya category of the Lefschetz fibration defined by a Brieskorn-Pham polynomial is equivalent to the triangulated category of singularities associated with the same polynomial together with a grading by an abelian group of rank one. Symplectic Picard-Lefschetz theory developed by Seidel is an essential ingredient of the proof.

## 1 Introduction

A polynomial $f \in \mathbb{C}[x_1, \dots, x_n]$ is said to be a *Brieskorn-Pham polynomial* if

$$f = f_{\boldsymbol{p}} = x_1^{p_1} + \cdots + x_n^{p_n}$$

for a sequence $\boldsymbol{p} = (p_1, \dots, p_n)$ of positive integers. A hypersurface singularity defined by a Brieskorn-Pham polynomial is called a *Brieskorn-Pham singularity*. This class of singularities includes a part of simple singularities, simple elliptic singularities, Arnold's exceptional unimodal singularities, and many more.

An important invariant of a hypersurface singularity is the *Milnor lattice*, which is the homology group of the Milnor fiber equipped with the intersection form. More recently, Seidel [21] introduced the *Fukaya category of a Lefschetz fibration*, which is a categorification of the Milnor lattice in the sense that the Grothendieck group equipped with the symmetrized Euler form is naturally isomorphic to the Milnor lattice.

Although the Milnor lattice of a singularity is difficult to compute in general, the Milnor lattice of a Brieskorn-Pham singularity allows the following description: The Milnor lattice of an $A_{p-1}$-singularity, defined by $f_p = x^p$ for an integer $p$ greater than one, is a free Abelian group generated by $C_i$ for $i = 1, \dots, p-1$ with the intersection form given by

$$(C_i, C_j) = \begin{cases} 2 & i = j, \\ -1 & |i - j| = 1, \\ 0 & \text{otherwise.} \end{cases}$$

Note that this is isomorphic to the root lattice of type $A_{p-1}$. Let

$$I_{\boldsymbol{p}} = \{(i_1, \dots, i_n) \in \mathbb{N}^n \mid 1 \le i_k \le p_k - 1, \ k = 1, \dots, n\},$$

be a finite set equipped with the lexicographic order

$$(i_1, \dots, i_n) < (j_1, \dots, j_n) \quad \text{if } i_k = j_k \text{ for } k < \ell \text{ and } i_\ell < j_\ell \text{ for some } \ell \in \{1, \dots, n\}.$$

Then it follows from a theorem of Sebastiani and Thom [17] that the Milnor lattice of $f_{\boldsymbol{p}} = x_1^{p_1} + \cdots + x_n^{p_n}$ is the *tensor product* of the Milnor lattices for $f_{p_k}$, so that there is a distinguished basis $(C_{\boldsymbol{i}})_{\boldsymbol{i} \in I_{\boldsymbol{p}}}$ of vanishing cycles satisfying

$$(C_{\boldsymbol{i}}, C_{\boldsymbol{j}}) = \begin{cases} \prod_{k=1}^n (C_{i_k}, C_{j_k}) & \text{if } i_k \le j_k \text{ for } k = 1, \dots, n, \\ 0 & \text{otherwise} \end{cases} \tag{1.1}$$

for $\boldsymbol{i} < \boldsymbol{j}$. The intersection form is determined by (1.1) together with the anti-symmetry

$$(C_{\boldsymbol{i}}, C_{\boldsymbol{j}}) = -(C_{\boldsymbol{j}}, C_{\boldsymbol{i}})$$

if $n$ is even, and the symmetry

$$(C_{\boldsymbol{i}}, C_{\boldsymbol{j}}) = (C_{\boldsymbol{j}}, C_{\boldsymbol{i}})$$

and $(C_{\boldsymbol{i}}, C_{\boldsymbol{i}}) = 2$ if $n$ is odd.

The main result in this paper is a categorification of the above description of the Milnor lattice. For an integer $p$ greater than one, let $\mathfrak{A}_p$ be the differential graded category whose set of objects is

$$\mathfrak{Ob}(\mathfrak{A}_p) = (C_1, \dots, C_p),$$

and whose spaces of morphisms are

$$\hom(C_i, C_j) = \begin{cases} \mathbb{C} \cdot \mathrm{id}_{C_i} & \text{if } i = j, \\ \mathbb{C}[-1] & \text{if } i = j - 1, \\ 0 & \text{otherwise,} \end{cases}$$

with the trivial differential. The tensor product of differential graded categories $\mathcal{A}$ and $\mathcal{B}$ is defined by

$$\mathfrak{Ob}(\mathcal{A} \otimes \mathcal{B}) = \mathfrak{Ob}(\mathcal{A}) \times \mathfrak{Ob}(\mathcal{B})$$

and

$$\mathrm{Hom}_{\mathcal{A} \otimes \mathcal{B}}(A_1 \times B_1, A_2 \times B_2) = \mathrm{Hom}_{\mathcal{A}}(A_1, A_2) \otimes \mathrm{Hom}_{\mathcal{B}}(B_1, B_2)$$

together with the differential determined by the Leibniz rule. Now consider the polynomial

$$W_{\boldsymbol{p}}(x_1, \dots, x_n) = f_{\boldsymbol{p}}(x_1, \dots, x_n) + (\text{lower order terms})$$

obtained by Morsifying $f_{\boldsymbol{p}}$ and let $\mathfrak{Fuk}\, W_{\boldsymbol{p}}$ be the Fukaya category of $W_{\boldsymbol{p}} : \mathbb{C}^n \to \mathbb{C}$ considered as an exact symplectic Lefschetz fibration with respect to the Euclidean Kähler structure on $\mathbb{C}^n$.

**Theorem 1.1.** *For any sequence $\boldsymbol{p} = (p_1, \dots, p_n)$ of positive integers, one has a quasi-equivalence*

$$\mathfrak{Fuk}\, W_{\boldsymbol{p}} \cong \mathfrak{A}_{p_1 - 1} \otimes \cdots \otimes \mathfrak{A}_{p_n - 1}$$

*of $A_\infty$-categories.*

An essential ingredient of the proof is the symplectic Picard-Lefschetz theory developed by Seidel [21], which provides an inductive tool to compute the Fukaya category in a combinatorial way.

Another important category that one can associate with a singularity is the *stabilized derived category*, introduced by Buchweitz [1] as the quotient category

$$D_{\mathrm{Sg}}(A) = D^b(\operatorname{mod} A)/D^{\mathrm{perf}}(\operatorname{mod} A)$$

of the bounded derived category $D^b(\operatorname{mod} A)$ of finitely-generated $A$-modules by its full subcategory $D^{\mathrm{perf}}(\operatorname{mod} A)$ consisting of perfect complexes. Here $A$ is the coordinate ring of the singularity, and a complex of $A$-module is said to be *perfect* if it is quasi-isomorphic to a bounded complex of projective modules. The motivation for this category comes from *matrix factorizations*, introduced by Eisenbud [4] to study maximal Cohen-Macaulay modules on a hypersurface. The same category is studied by Orlov [13] under the name '*triangulated category of singularities*'.

If the ring $A$ is graded by an abelian group, then there is a graded version of the stabilized derived category, defined as the quotient category

$$D^{\mathrm{gr}}_{\mathrm{Sg}}(A) = D^b(\operatorname{gr} A)/D^{\mathrm{perf}}(\operatorname{gr} A)$$

of the bounded derived category $D^b(\operatorname{gr} A)$ of finitely-generated graded $A$-modules by its full subcategory $D^{\mathrm{perf}}(\operatorname{gr} A)$ consisting of bounded complexes of projectives.

In the case of a Brieskorn-Pham singularity, we equip the coordinate ring

$$A_{\boldsymbol{p}} = \mathbb{C}[x_1, \dots, x_n]/(f_{\boldsymbol{p}}),$$

with the grading given by the abelian group $L(\boldsymbol{p})$ of rank one generated by $n+1$ elements $\vec{x}_1, \dots, \vec{x}_n, \vec{c}$ with relations

$$p_1\vec{x}_1 = p_2\vec{x}_2 = \cdots = p_n\vec{x}_n = \vec{c}.$$

**Theorem 1.2.** *For any sequence $\boldsymbol{p} = (p_1, \dots, p_n)$ of positive integers, one has an equivalence*

$$D^{\mathrm{gr}}_{\mathrm{Sg}}(A_{\boldsymbol{p}}) \cong D^b(\mathfrak{A}_{p_1-1} \otimes \cdots \otimes \mathfrak{A}_{p_n-1})$$

*of triangulated categories.*

By combining Theorem 1.1 and Theorem 1.2, one obtains homological mirror symmetry for Brieskorn-Pham singularities:

**Theorem 1.3.** *For any sequence $\boldsymbol{p} = (p_1, \dots, p_n)$ of positive integers, one has an equivalence*

$$D^b \operatorname{\mathfrak{Fuk}} W_{\boldsymbol{p}} \cong D^{\mathrm{gr}}_{\mathrm{Sg}}(A_{\boldsymbol{p}})$$

*of triangulated categories.*

The special case of $n = 2$ in Theorem 1.3 is proved in [25]. Homological mirror symmetry is proposed by Kontsevich [10] for Calabi-Yau manifolds, and later generalized to more general classes of manifolds [11, 8, 18]. Ebeling and Takahashi [3, 23, 22] discuss the relation between mirror symmetry for singularities and Saito's duality for regular systems of weights [16]. Stabilized derived categories of singularities associated with regular systems of weights whose smallest exponents are $\pm 1$ are studied by Kajiura, Saito and Takahashi [24, 6, 7]. Okada [12] also discusses homological mirror symmetry for Brieskorn-Pham singularities.

If $\boldsymbol{p}$ satisfies a condition described below, then one can relate the stabilized derived category $D^{\mathrm{gr}}_{\mathrm{Sg}}(A_{\boldsymbol{p}})$ with the derived category $D^b \operatorname{coh} Y_{\boldsymbol{p}}$ of coherent sheaves on a stack $Y_{\boldsymbol{p}}$ defined as follows: Let $\ell$ be the least common multiple of $(p_1, \dots, p_n)$ and equip $A(\boldsymbol{p})$ with a $\mathbb{Z}$-grading given by

$$\deg x_i = a_i = \frac{\ell}{p_i}, \qquad i = 1, \dots, n.$$

Then $X_{\boldsymbol{p}} = \operatorname{Proj} A_{\boldsymbol{p}}$ is a hypersurface of degree $\ell$ in the weighted projective space $\mathbb{P}(a_1, \dots, a_n)$. Put

$$K = \{(\alpha_1, \dots, \alpha_n) \in (\mathbb{C}^\times)^n \mid \alpha_1^{p_1} = \cdots = \alpha^{p_n}\}$$

and define a homomorphism

$$\phi : \mathbb{C}^\times \to K$$

by

$$\phi(\alpha) = (\alpha^{a_1}, \dots, \alpha^{a_n}).$$

Then the cokernel

$$G_{\boldsymbol{p}} = \operatorname{coker} \phi$$

of $\phi$ is a finite abelian group acting on $X_{\boldsymbol{p}}$, and let

$$Y_{\boldsymbol{p}} = [X_{\boldsymbol{p}} / G_{\boldsymbol{p}}]$$

be the quotient stack with respect to this action. The adaptation of the Calabi-Yau/Landau-Ginzburg correspondence proved by Orlov [14, Theorem 2.5] to the $L(\boldsymbol{p})$-graded situation gives the following:

**Theorem 1.4.** *If a sequence $\boldsymbol{p} = (p_1, \dots, p_n)$ satisfies*

$$\frac{1}{p_1} + \cdots + \frac{1}{p_n} = 1,$$

*then one has an equivalence*

$$D^{\mathrm{gr}}_{\mathrm{Sg}}(A_{\boldsymbol{p}}) \cong D^b \operatorname{coh} Y_{\boldsymbol{p}}$$

*of triangulated categories.*

By combining Theorem 1.3 and Theorem 1.4, one obtains an equivalence

$$D^b \mathfrak{Fuk}\, W_{\boldsymbol{p}} \cong D^b \operatorname{coh}^G Y$$

between the derived Fukaya category of the Lefschetz fibration $W_{\boldsymbol{p}}$ and the derived category of coherent sheaves on the stack $Y_{\boldsymbol{p}}$.

The organization of this paper is as follows: In Section 2, we give a description of the Fukaya category of $f(x) + u^k$ in terms of the Fukaya category of $f^{-1}(0)$. This is based on symplectic Picard-Lefschetz theory developed by Seidel [21], and the case when $k = 2$ is discussed in [19]. In Section 3, we prove Theorem 1.1 by induction on $n$. The proof of Theorem 1.2 is given in Section 4.

*Acknowledgment*: M. F. is supported by Grant-in-Aid for Young Scientists (No.19.8083). K. U. is supported by Grant-in-Aid for Young Scientists (No.18840029). This work has been done while K. U. is visiting the University of Oxford, and he thanks the Mathematical Institute for hospitality and Engineering and Physical Sciences Research Council for financial support.

# 2 The Fukaya category of $f(x) + u^k$

Let $f : \mathbb{C}^n \to \mathbb{C}$ be a polynomial in $n$ variables. Assume that

- $f$ is *tame*, in the sense that the gradient $\|\nabla f\|$ is bounded from below outside of a compact set by a positive number, and
- $f$ has non-degenerate critical points with distinct critical values.

Then $f$ gives an *exact symplectic Lefshetz fibration* [21, Section (15d)] with respect to the standard Euclidean Kähler structure on $\mathbb{C}^n$. Assume for simplicity that the set of critical values of $f$ is the set of $m$-th roots of unity, and let $(\gamma)_{i=1}^m$ be the distinguished set of vanishing paths chosen as straight line segments from the origin. Let $\mathcal{B}$ denote the Fukaya category of $f^{-1}(0)$ whose objects are vanishing cycles $C_i \subset f^{-1}(0)$ along $\gamma_i$ for $i = 1, \dots, m$. It is an $A_\infty$-category whose spaces of morphisms are given by Lagrangian intersection Floer complexes, and $A_\infty$-operations are defined by counting virtual numbers of stable maps of genus zero with Lagrangian boundary conditions [5]. The directed subcategory of $\mathcal{B}$ with respect to the order $(C_1, \dots, C_m)$ will be denoted by $\mathcal{A}$. Although $\mathcal{A}$ depends on the choice of a distinguished basis of vanishing cycles, the derived category $D^b\mathcal{A}$ is independent of this choice [20, 21] and gives an invariant of the Lefschetz fibration.

Let

$$g(u) = u^k - \epsilon u$$

be a perturbation of $u^k$ and consider the *Lefschetz bifibration* [21, Section (15e)]

$$\Psi = \psi \circ \varpi$$

$$\mathbb{C}^{n+1} \xrightarrow{\ \varpi\ } \mathbb{C}^2 \xrightarrow{\ \psi\ } \mathbb{C}$$

where

$$\varpi(x, u) = (f(x) + g(u), u)$$

and

$$\psi(y_1, y_2) = y_1.$$

We write the critical points of $f$ and $g$ as

$$\operatorname{Crit} f = \{x_i\}_{i=1}^m$$

and

$$\operatorname{Crit} g = \{u \in \mathbb{C} \mid k u^{k-1} = \epsilon\} = \{u_j\}_{j=1}^{k-1},$$

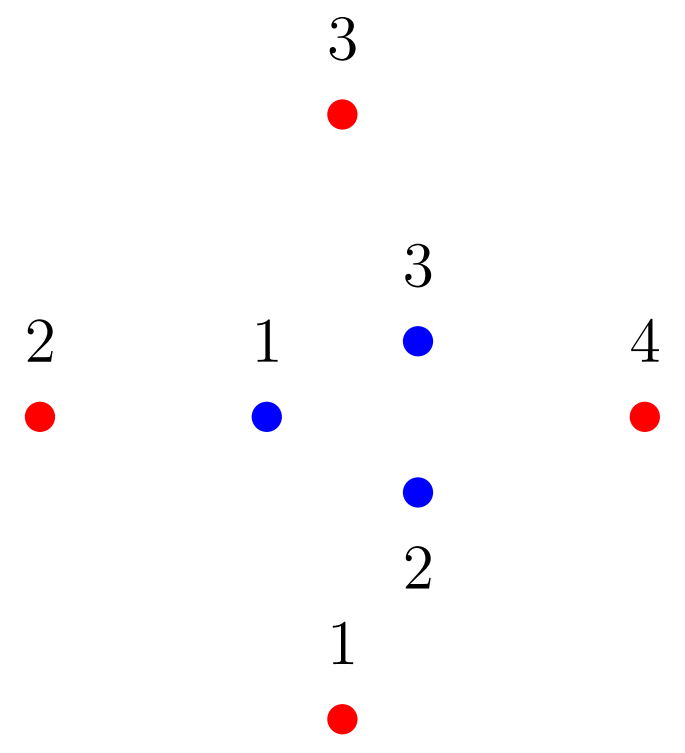

Figure 2.1: Critical values of $f$ (outside) and $g$ (inside)

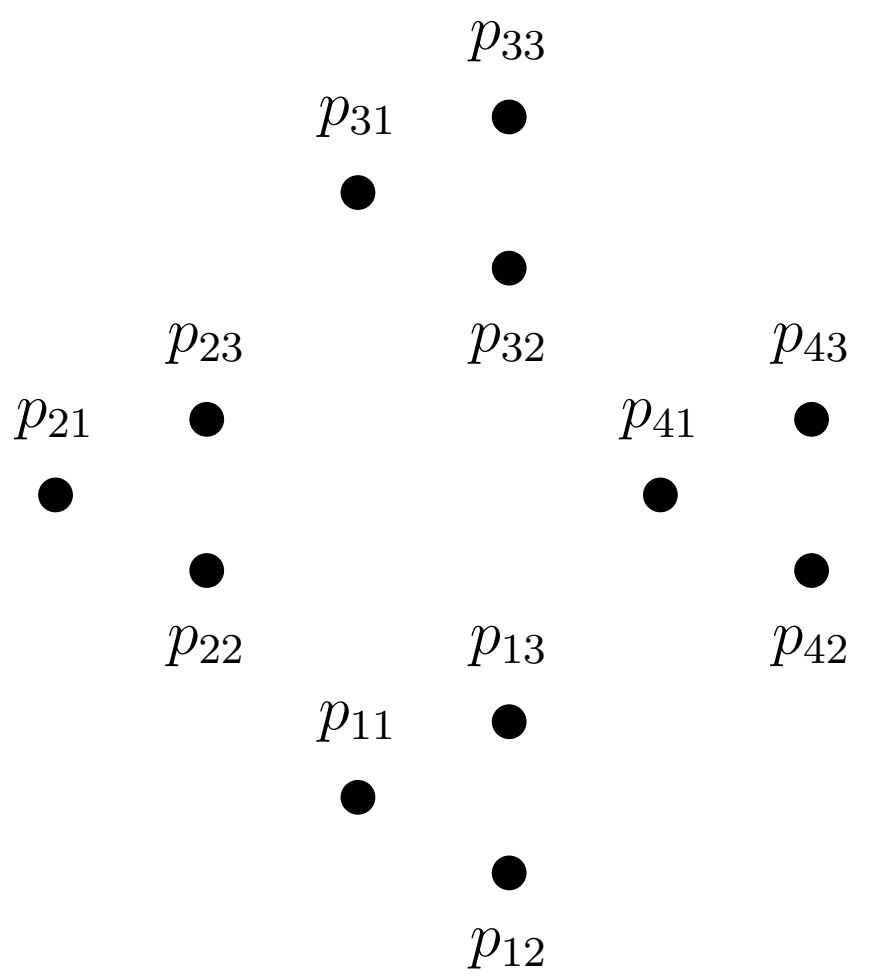

Figure 2.2: Critical values of $\Psi$

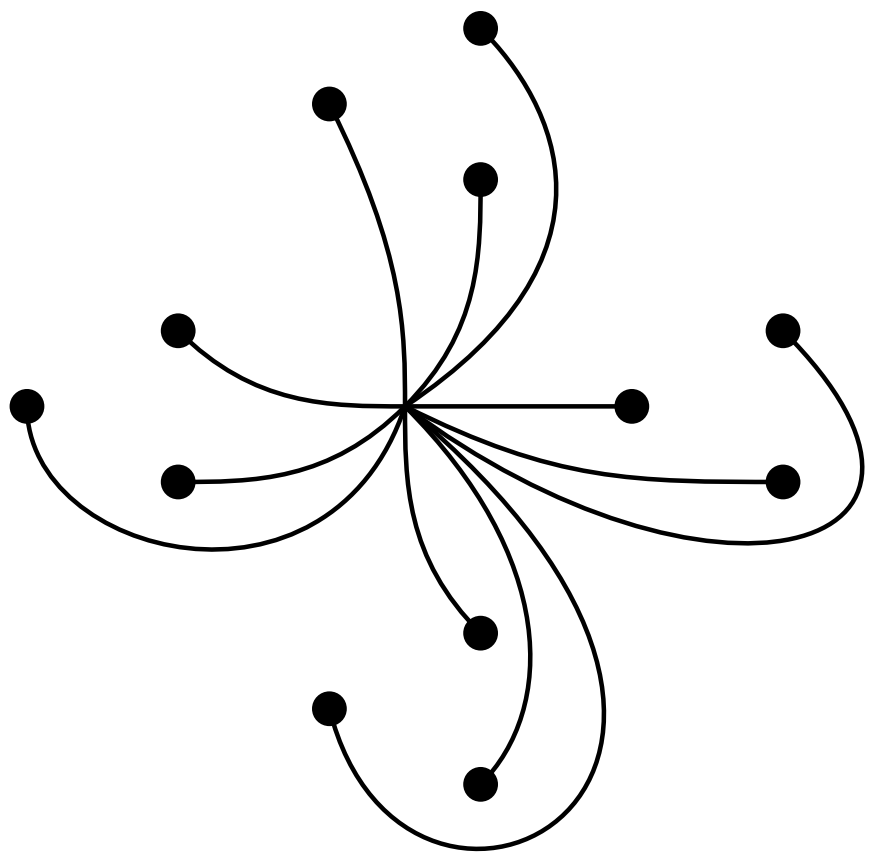

Figure 2.3: A distinguished set of vanishing paths for $\Psi$

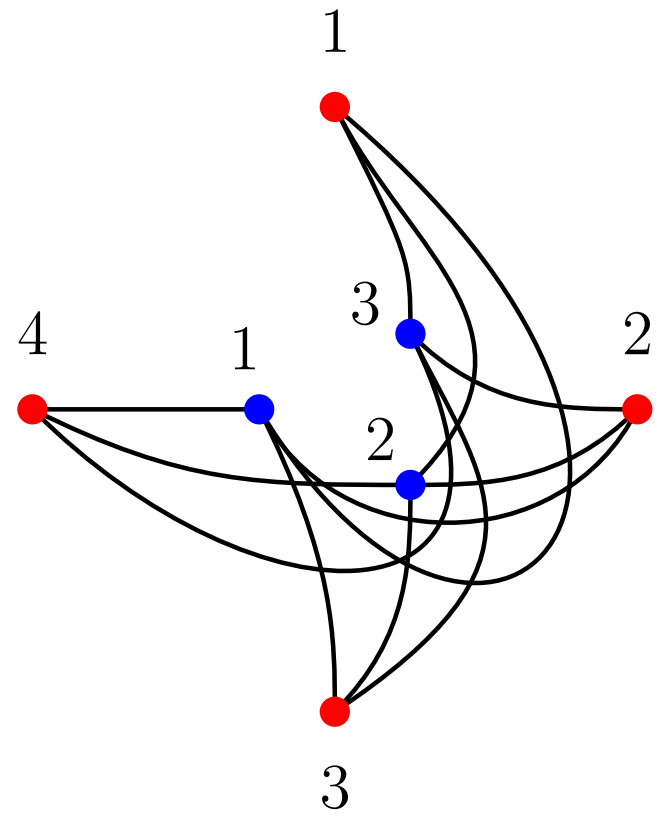

Figure 2.4: Trajectories of critical values of $t - f$ along vanishing paths

so that the set of critical points of $\Psi$ is given by

$$\operatorname{Crit}\Psi = \operatorname{Crit} f \times \operatorname{Crit} g = \{p_{ij} = (x_i, u_j)\}_{i,j}$$

with critical values

$$\Psi(p_{ij}) = f(x_i) + g(u_j).$$

Assume that $\epsilon$ is sufficiently small and order the critical values of $f$ clockwise and those of $g$ counterclockwise as shown in Figure 2.1. The critical values of $\Psi$ are shown in Figure 2.2, and we choose a distinguished set of vanishing paths $\gamma_{ij}$ from the origin to $p_{ij}$ as in Figure 2.3. These figures are for the case $m = k = 4$, and the general case is similar.

For general $t \in \mathbb{C}$, the map

$$\mathcal{E}_t \xrightarrow{\varpi_t} \mathcal{S}_t$$

from $\mathcal{E}_t = \Psi^{-1}(t)$ to $\mathcal{S}_t = \psi^{-1}(t)$ is a Lefschetz fibration. The vanishing cycle $C_{ij}$ of $\Psi$ along the path $\gamma_{ij}$ comes from a matching path $\mu_{ij}$, which is obtained as the trajectory of critical values of $\varpi_t$ along the path $\gamma_{ij}$.

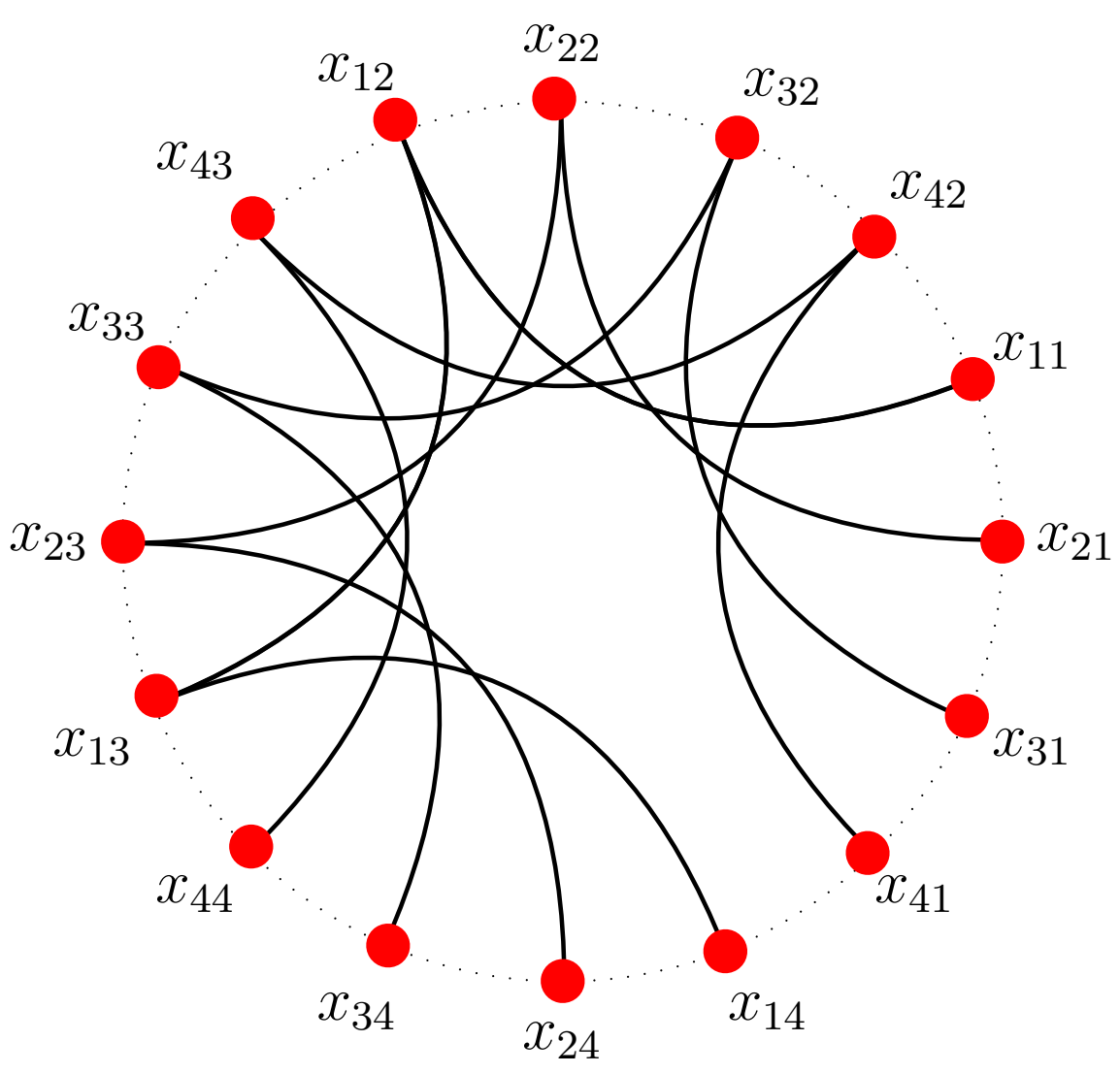

Figure 2.5: Matching paths on the $u$-plane

The fiber $\mathcal{S}_t$ of $\psi$ can be identified with the $u$-plane, and the image of the matching path by $g$ can be described as follows: Consider the map

$$\begin{array}{cccc} \pi_t: & \Psi^{-1}(t) & \to & \mathbb{C} \\ & \Downarrow & & \Downarrow \\ & (x,u) & \mapsto & s = t - f(x) = g(u). \end{array}$$

The fiber of $\pi_t$ over $s \in \mathbb{C}$ is given by

$$\pi_t^{-1}(s) = f^{-1}(t-s) \times g^{-1}(s).$$

The first factor becomes singular if $s$ is $t$ minus a critical value of $f$, and the second factor becomes singular if $s$ is a critical value of $g$. As one varies $t$ along vanishing paths in Figure 2.3, the critical values of minus $f$ move on the $s$-plane until one of them hits one of the critical values of $g$. The trajectory of $t$ minus a critical value of $f$ along $\gamma_{ij}$, starting from the $i$-th point outside and ending at the $j$-th point inside, is the image by $g$ of the matching path $\mu_{ij}$ on the $u$-plane corresponding to the vanishing path $\gamma_{ij}$. Figure 2.4 shows these trajectories, and Figure 2.5 shows the matching paths obtained as the inverse images of these trajectories by $g$.

For each critical point $x_i$ of $f$, there are $k$ critical values $g^{-1} \circ f(x_i) = \{x_{i1}, \dots, x_{ik}\}$ of $\varpi_0$, and the matching path $\mu_{ij}$ connects $x_{i,j+1}$ with $x_{ij}$. We write the straight line segment on the $u$-plane from the origin to $x_{ij}$ as $\delta_{ij}$. The fiber $\varpi_0^{-1}(0)$ can naturally be identified with $f^{-1}(0)$, so that the vanishing cycle $\Delta_{i,j}$ of $\varpi_0$ along $\delta_{ij}$ corresponds to $C_i$. This shows that the Fukaya category $\mathcal{B}_k$ of $\varpi_0^{-1}(0)$ consisting of $\Delta_{i,j}$ is given by

$$\hom_{\mathcal{B}_k}(\Delta_{i,j}, \Delta_{i',j'}) = \hom_{\mathcal{B}}(C_i, C_{i'})$$

with the natural $A_\infty$-structure inherited from $\mathcal{B}$.

Let $\mathcal{A}_k$ be the directed subcategory of $\mathcal{B}_k$ with respect to the order

$$(i, j) < (i', j') \qquad \text{if} \quad j > j' \quad \text{or} \quad j = j' \text{ and } i < i'$$

on the index set. Note that the $A_\infty$-category $\mathcal{A}_k$ depends not only on $\mathcal{A}$ but also on $\mathcal{B}$. Let further

$$S_{ij} = \operatorname{Cone}(\Delta_{i,j+1} \to \Delta_{i,j}), \qquad i = 1, \dots, m, \quad j = 1, \dots, k-1$$

be the object in $D^b\mathcal{A}_k$ which is the cone over the morphism $e_{i,j} : \Delta_{i,j+1} \to \Delta_{i,j}$ corresponding to $\mathrm{id}_{C_i}$ under the isomorphism

$$\hom_{\mathcal{A}_k}(\Delta_{i,j+1}, \Delta_{i,j}) \cong \hom_{\mathcal{B}}(C_i, C_i).$$

The following theorem gives a description of the Fukaya category of $\Psi^{-1}(0)$ in terms of the Fukaya category of $f^{-1}(0)$:

**Theorem 2.1** (Seidel [21, Proposition 18.21])**.** *If $n$ is greater than one, then the Fukaya category of $\Psi^{-1}(0)$ consisting of vanishing cycles $C_{ij}$ is quasi-equivalent to the full subcategory of $D^b\mathcal{A}_k$ consisting of $S_{ij}$.*

# 3 Inductive description of the Fukaya category

Let $\boldsymbol{p} = (p_1, \dots, p_n)$ be a sequence of natural numbers and $\widetilde{\boldsymbol{p}} = (p_1, \dots, p_n, k)$ be another sequence obtained by appending $p_{n+1} = k$ to $\boldsymbol{p}$. Let further $W_{\boldsymbol{p}}$ be a perturbation of the Brieskorn-Pham polynomial of degree $\boldsymbol{p}$ and $g(u)$ be a perturbation of $u^k$ as in Section 2, so that $W_{\widetilde{\boldsymbol{p}}} = W_{\boldsymbol{p}} + g$ is a perturbation of the Brieskorn-Pham polynomial of degree $\widetilde{\boldsymbol{p}}$. The directed Fukaya category $\mathfrak{Fuk}\, W_{\boldsymbol{p}}$ consisting of a distinguished basis of vanishing cycles in $W_{\boldsymbol{p}}^{-1}(0)$ will be denoted by $\mathcal{A}$.

Assume that Theorem 1.1 holds for $W_{\boldsymbol{p}}$ so that one has a quasi-equivalence

$$\mathcal{A} = \mathfrak{A}_{p_1-1} \otimes \cdots \otimes \mathfrak{A}_{p_n-1}$$

of $A_\infty$-categories. Theorem 2.1 shows that $\mathfrak{Fuk}\, W_{\widetilde{\boldsymbol{p}}}$ is quasi-equivalent to the directed subcategory of $D^b\mathcal{A}_k$ consisting of

$$C_{\boldsymbol{i},j} = \operatorname{Cone}(\Delta_{\boldsymbol{i},j+1} \to \Delta_{\boldsymbol{i},j}), \qquad (\boldsymbol{i}, j) \in I_{\widetilde{\boldsymbol{p}}}$$

with respect to the order

$$(\boldsymbol{i}, j) < (\boldsymbol{i}', j') \qquad \text{if} \quad \boldsymbol{i} < \boldsymbol{i}' \quad \text{or} \quad \boldsymbol{i} = \boldsymbol{i}' \text{ and } j < j'.$$

Note that one has

$$\hom_{D^b\mathcal{A}_k}(C_{\boldsymbol{i},j}, C_{\boldsymbol{i}',j'}) = \left\{ \begin{array}{ccc} \hom_{\mathcal{A}_k}(\Delta_{\boldsymbol{i},j}, \Delta_{\boldsymbol{i}',j'+1}) & \xrightarrow{\bullet \circ e_{\boldsymbol{i},j}} & \hom_{\mathcal{A}_k}(\Delta_{\boldsymbol{i},j+1}, \Delta_{\boldsymbol{i}',j'+1}) \\ {\scriptstyle e_{\boldsymbol{i}',j'} \circ \bullet} \downarrow & & {\scriptstyle e_{\boldsymbol{i}',j'} \circ \bullet} \downarrow \\ \hom_{\mathcal{A}_k}(\Delta_{\boldsymbol{i},j}, \Delta_{\boldsymbol{i}',j'}) & \xrightarrow{\bullet \circ e_{\boldsymbol{i},j}} & \hom_{\mathcal{A}_k}(\Delta_{\boldsymbol{i},j+1}, \Delta_{\boldsymbol{i}',j'}) \end{array} \right\},$$

where the right hand side denotes the total complex of the double complex. If $j < j' - 1$, then the right hand side is trivial;

$$\left\{ \begin{array}{ccc} 0 & \longrightarrow & 0 \\ \big\downarrow & & \big\downarrow \\ 0 & \longrightarrow & 0 \end{array} \right\} \simeq 0.$$

If $j = j' - 1$, then the right hand side is given by

$$\left\{ \begin{array}{ccc} 0 & \longrightarrow & 0 \\ \big\downarrow & & \big\downarrow \\ 0 & \longrightarrow & \hom_{\mathcal{A}}(C_{\boldsymbol{i}}, C_{\boldsymbol{i}'}) \end{array} \right\} \simeq \hom_{\mathcal{A}}(C_{\boldsymbol{i}}, C_{\boldsymbol{i}'})[-1].$$

Natural representatives of a basis of the cohomology group of this complex are given by

$$\begin{array}{ccccc} & & \Delta_{\boldsymbol{i},j+1} & \longrightarrow & \Delta_{\boldsymbol{i},j} \\ & & {\scriptstyle x_{\boldsymbol{i},\boldsymbol{i}'}}\big\downarrow & & \\ \Delta_{\boldsymbol{i}',j+2} & \longrightarrow & \Delta_{\boldsymbol{i}',j+1} & & \end{array}$$

where $x_{\boldsymbol{i},\boldsymbol{i}'}$ runs over a basis of $\hom_{\mathcal{A}}(C_{\boldsymbol{i}}, C_{\boldsymbol{i}'})$. If $j = j'$ and $\boldsymbol{i} < \boldsymbol{i}'$, then the right hand side is given by

$$\left\{ \begin{array}{ccc} 0 & \longrightarrow & \hom_{\mathcal{A}}(C_{\boldsymbol{i}}, C_{\boldsymbol{i}'}) \\ \big\downarrow & & {\scriptstyle \mathrm{id}}\big\downarrow \\ \hom_{\mathcal{A}}(C_{\boldsymbol{i}}, C_{\boldsymbol{i}'}) & \xrightarrow{\mathrm{id}} & \hom_{\mathcal{A}}(C_{\boldsymbol{i}}, C_{\boldsymbol{i}'}) \end{array} \right\} \simeq \hom_{\mathcal{A}}(C_{\boldsymbol{i}}, C_{\boldsymbol{i}'}),$$

whose cohomologies are spanned by

$$\begin{array}{ccc} \Delta_{\boldsymbol{i},j+1} & \longrightarrow & \Delta_{\boldsymbol{i},j} \\ {\scriptstyle x_{\boldsymbol{i},\boldsymbol{i}'}}\big\downarrow & & {\scriptstyle x_{\boldsymbol{i},\boldsymbol{i}'}}\big\downarrow \\ \Delta_{\boldsymbol{i}',j+1} & \longrightarrow & \Delta_{\boldsymbol{i}',j}. \end{array}$$

If $j > j'$ and $\boldsymbol{i} < \boldsymbol{i}'$, then the right hand side is acyclic;

$$\left\{ \begin{array}{ccc} \hom_{\mathcal{A}}(C_{\boldsymbol{i}}, C_{\boldsymbol{i}'}) & \xrightarrow{\mathrm{id}} & \hom_{\mathcal{A}}(C_{\boldsymbol{i}}, C_{\boldsymbol{i}'}) \\ {\scriptstyle \mathrm{id}}\big\downarrow & & {\scriptstyle \mathrm{id}}\big\downarrow \\ \hom_{\mathcal{A}}(C_{\boldsymbol{i}}, C_{\boldsymbol{i}'}) & \xrightarrow{\mathrm{id}} & \hom_{\mathcal{A}}(C_{\boldsymbol{i}}, C_{\boldsymbol{i}'}) \end{array} \right\} \simeq 0.$$

It is straightforward to compute the compositions among the above basis to show that the cohomology category of $\mathfrak{Fuk}\, W_{\tilde{\boldsymbol{p}}}$ is equivalent to $\mathfrak{A}_{p_1-1} \otimes \cdots \otimes \mathfrak{A}_{p_n-1} \otimes \mathfrak{A}_{k-1}$ as a graded category. Higher $A_\infty$-operations on $\mathfrak{Fuk}\, W_{\tilde{\boldsymbol{p}}}$ vanish for degree reasons, and one obtains a quasi-equivalence

$$\mathfrak{Fuk}\, W_{\tilde{\boldsymbol{p}}} \cong \mathfrak{A}_{p_1-1} \otimes \cdots \otimes \mathfrak{A}_{p_n-1} \otimes \mathfrak{A}_{k-1},$$

of $A_\infty$-categories. It is easy to check the $n = 2$ case using Figure 2.5 since only triangles contribute because of the directedness, and Theorem 1.1 is proved.

# 4 The triangulated category of singularities

Let $V$ be a vector space and $W$ be an element of the polynomial ring $\mathbb{C}[V]$. Assume that $W$ has an isolated critical point at the origin and there is a one-form $\gamma$ such that

$$W = \gamma(\eta),$$

where $\eta = \sum_i x_i \partial_i$ is the Euler vector field. The following is a variation of [18, Lemma 12.3]:

**Lemma 4.1.** *Let $\mathcal{P}$ be the structure sheaf of the origin as an $\mathcal{O}_X$-module, where $X$ is the zero locus of $W$. Then the chain complex $C = (C^i, \delta^i)$ given by*

$$C^i = \begin{cases} \bigoplus_{j=0}^{\lfloor -i/2 \rfloor} \Omega_V^{-i-2j}|_X & i \leq 0, \\ 0 & i > 0, \end{cases}$$
$$\delta^i = \iota_\eta + \gamma \wedge \cdot,$$

*is an $\mathcal{O}_X$-free resolution of $\mathcal{P}$.*

*Proof.* Let $\mathcal{Q}$ be the structure sheaf of the origin as an $\mathcal{O}_V$-module. Then the derived pull-back of $\mathcal{Q}$ by the inclusion $i : X \hookrightarrow V$ is the direct sum of $\mathcal{P}$ and $\mathcal{P}[1]$:

$$i^* \mathcal{Q} \cong \mathcal{P} \oplus \mathcal{P}[1].$$

Since the Koszul complex

$$K = \left\{ 0 \to \Omega_V^n \xrightarrow{\iota_\eta} \Omega_V^{n-1} \xrightarrow{\iota_\eta} \cdots \xrightarrow{\iota_\eta} \Omega_V^0 \to 0 \right\}$$

is a $\mathcal{O}_V$-free resolution of $\mathcal{Q}$, this shows that its restriction

$$K|_X = \left\{ 0 \to \Omega_V^n|_X \xrightarrow{\iota_\eta} \Omega_V^{n-1}|_X \xrightarrow{\iota_\eta} \cdots \xrightarrow{\iota_\eta} \Omega_V^0|_X \to 0 \right\}$$

to $X$ is also isomorphic to the direct sum of $\mathcal{P}$ and $\mathcal{P}[1]$. Now consider the chain map

$$\begin{array}{ccccccccccccccc} 0 & \to & \Omega_V^n|_X & \to & \Omega_V^{n-1}|_X & \to & \cdots & \to & \Omega_V^1|_X & \to & \Omega_V^0|_X & \to & 0 & \to & 0 \\ & & \downarrow & & \downarrow & & & & \downarrow & & \downarrow & & \downarrow & & \\ 0 & \to & 0 & \to & \Omega_V^n|_X & \to & \cdots & \to & \Omega_V^2|_X & \to & \Omega_V^1|_X & \to & \Omega_V^0|_X & \to & 0 \end{array}$$

from $K|_X[1]$ to $K|_X$ where vertical arrows are given by $\gamma \wedge \cdot$. Since this induces the identity map on the $(-1)$-th cohomology group, which is $\mathcal{P}$ for both $K|_X[1]$ and $K|_X$, the mapping cone for this map is isomorphic to $\mathcal{P} \oplus \mathcal{P}[2]$;

$$\left\{ K|_X[1] \xrightarrow{\gamma \wedge \cdot} K|_X \right\} \cong \mathcal{P} \oplus \mathcal{P}[2].$$

By iterating this process, one obtains

$$\left\{ K|_X[i] \xrightarrow{\gamma \wedge \cdot} K|_X[i-1] \xrightarrow{\gamma \wedge \cdot} \cdots \xrightarrow{\gamma \wedge \cdot} K|_X \right\} \cong \mathcal{P} \oplus \mathcal{P}[i+1].$$

Now the lemma follows by taking $i$ to infinity. $\square$

Now we prove Theorem 1.2 along the lines of [25, Theorem 5]. Fix any weight $\boldsymbol{p} = (p_1, \dots, p_n)$ and put $A = A(\boldsymbol{p})$ and $L = L(\boldsymbol{p})$. We will find a full triangulated subcategory $\mathcal{T}$ of $D^b(\operatorname{gr} A)$ equivalent to $D^{\mathrm{gr}}_{\mathrm{Sg}}(A)$ such that $(k(\vec{n}))_{\vec{n} \in I}$ is a full exceptional collection in $\mathcal{T}$, where $k = A/(x_1, \dots, x_n)$ and

$$I = \{a_1\vec{x}_1 + \cdots + a_n\vec{x}_n \in L \mid -p_1 + 2 \le a_1 \le 0, \dots, -p_n + 2 \le a_n \le 0\}.$$

Let $L_+$ be the subset of $L$ defined by

$$L_+ = \{-(n-1)\vec{c} + a_1\vec{x}_1 + \cdots + a_n\vec{x}_n \mid a_i \ge 1, \quad i = 1, \dots, n\},$$

and $L_-$ be the complement $L \setminus L_+$. Let further $\mathcal{S}_-$ and $\mathcal{P}_+$ be the full triangulated subcategories of $D^b(\operatorname{gr} A)$ generated by $k(\vec{n})$ for $\vec{n} \in L_+$ and $A(\vec{m})$ for $\vec{m} \in L_-$ respectively. Then $\mathcal{S}_-$ and $\mathcal{P}_+$ are left admissible in $D^b(\operatorname{gr} A)$, and since $\mathcal{P}_+ \subset {}^\perp\mathcal{S}_-$, one has a weak semiorthogonal decomposition

$$D^b(\operatorname{gr} A) = \langle \mathcal{S}_-, \mathcal{P}_+, \mathcal{T} \rangle$$

such that $\mathcal{T} \cong D^{\mathrm{gr}}_{\mathrm{Sg}}(A)$. One can see that $k(\vec{m})$ for $\vec{m} \in I$ belongs to $\mathcal{T}$, since

$$\mathbb{R}\operatorname{Hom}(k(\vec{m}), k(\vec{n})) = 0$$

if $\vec{m} \notin \vec{n} + \mathbb{N}\vec{x}_1 + \cdots + \mathbb{N}\vec{x}_n$, and

$$\mathbb{R}\operatorname{Hom}(k(\vec{m}), A(\vec{n})) = 0$$

if $\vec{m} \ne -\vec{c} + \vec{x}_1 + \cdots + \vec{x}_n + \vec{n}$. The $\mathbb{R}\operatorname{Hom}$'s between them can be calculated by the free resolution obtained in Lemma 4.1 to be

$$\begin{aligned} &\mathbb{R}\operatorname{Hom}(k(a_1\vec{x}_1 + \cdots + a_n\vec{x}_n), k(b_1\vec{x}_1 + \cdots + b_n\vec{x}_n)) \\ &\quad = \hom_{\mathcal{A}_{p_1-1}}(C_{-a_1+1}, C_{-b_1+1}) \otimes \cdots \otimes \hom_{\mathcal{A}_{p_n-1}}(C_{-a_n+1}, C_{-b_n+1}) \end{aligned}$$

Hence $(k(\vec{n}))_{\vec{n} \in I}$ is an exceptional collection. It is straightforward to read off the structure of the Yoneda products from the above resolution to show that the the full subcategory of $\mathcal{T}$ consisting of $(k(\vec{n}))_{\vec{n} \in I}$ is isomorphic as a graded category to $\mathfrak{A}_{p_1-1} \otimes \cdots \otimes \mathfrak{A}_{p_n-1}$. Moreover, $\mathcal{T}$ has a differential graded enhancement induced from that of $D^b(\operatorname{gr} A)$, which is formal for degree reasons. This shows that $\mathfrak{A}_{p_1-1} \otimes \cdots \otimes \mathfrak{A}_{p_n-1}$ is equivalent to the full triangulated subcategory of $D^{\mathrm{gr}}_{\mathrm{Sg}}(A)$ generated by the above exceptional collection.

To prove that the image of $(k(\vec{n}))_{\vec{n} \in I}$ in $D^{\mathrm{gr}}_{\mathrm{Sg}}(A)$ is full, we use the following:

**Lemma 4.2.** *The module $k(\vec{m})$ for any $\vec{m} \in L$ can be obtained from $(k(\vec{n}))_{\vec{n} \in I}$ by taking cones up to perfect complexes.*

*Proof.* First note that the exact sequences

$$\begin{aligned} 0 \to k(-\vec{x}_1) \to \quad & k[x_1]/(x_1^2) \quad \to k \to 0, \\ 0 \to k(-2\vec{x}_1) \to \quad & k[x_1]/(x_1^3) \quad \to k[x_1]/(x_1^2) \to 0, \\ & \vdots \\ 0 \to k(-(p_1-1)\vec{x}_1) \to \quad & k[x_1]/(x_1^{p_1}) \quad \to k[x_1]/(x_1^{p_1-1}) \to 0 \end{aligned}$$

of $A$-modules show that $k(-(p_1-1)\vec{x}_1)$ can be obtained from $k, k(-\vec{x}_1), \ldots, k(-(p_1-2)\vec{x}_1)$ by taking cones up to the perfect module $k[x_1]/(x_1^{p_1}) \cong A/(x_2, \ldots, x_n)$. Then by shifting the degrees, one can see that for any $\vec{n} \in L$, $k(\vec{n})$ can be obtained from either $k(\vec{n}-\vec{x}_1), k(\vec{n}-2\vec{x}_1), \ldots, k(\vec{n}-(p_1-1)\vec{x}_1)$ or $k(\vec{n}+\vec{x}_1), k(\vec{n}+2\vec{x}_1), \ldots, k(\vec{n}+(p_1-1)\vec{x}_1)$ by taking cones up to perfect complexes. The same is true for $\vec{x}_i$ for $i = 2, \ldots, n$, and the lemma follows. □

Now one can use [2, Corollary 4.3], [7, Theorem 4.5], [9, Proposition A.2], [15, Proposition 2.7] or [25, Section 4] to conclude that $(k(\vec{n}))_{\vec{n}\in I}$ is full.

Masahiro Futaki

Graduate School of Mathematical Sciences, The University of Tokyo, 3-8-1 Komaba Meguro-ku Tokyo 153-8914, Japan

*e-mail address* : futaki@ms.u-tokyo.ac.jp

Kazushi Ueda
Department of Mathematics, Graduate School of Science, Osaka University, Machikaneyama 1-1, Toyonaka, Osaka, 560-0043, Japan.
*e-mail address* : kazushi@math.sci.osaka-u.ac.jp